\documentclass[openacc]{rsproca_new}

\usepackage{microtype,marvosym} 

\usepackage[numbers]{natbib}

\usepackage{tikz,pgfplots,todonotes}
\pgfplotsset{compat=newest}
\pgfplotsset{plot coordinates/math parser=false}       

\usepackage{amsmath,amssymb,graphicx,MnSymbol}

\definecolor{lred}{RGB}{200,0,0}
\definecolor{dred}{RGB}{130,0,0} \definecolor{dblu}{RGB}{0,0,130}
\definecolor{dgre}{RGB}{0,130,0} \definecolor{dgra}{RGB}{50,50,50}
\definecolor{mgra}{RGB}{100,100,100}
\definecolor{lgra}{RGB}{220,220,220}
\definecolor{MPG}{RGB}{000,125,122}
\definecolor{ora}{HTML}{FF9933}

\newcommand{\mpg}[1]{{\color{MPG} #1}}  
  
\newcommand{\dre}[1]{{\color{dred} #1}}

\usepackage{booktabs}

\newcommand{\g}{\mid}

\newcommand{\cov}{\operatorname{cov}}

\renewcommand{\Re}{\mathbb{R}}

\newcommand{\N}{\mathcal{N}}

\newcommand{\Trans}{^{\intercal}}

\newcommand{\model}{\mathcal{M}}
\newcommand{\data}{\mathcal{D}}

\newcommand{\qq}{\qquad}

\renewcommand{\vec}{\boldsymbol} 
 
\newcommand{\vect}[1]{\overrightarrow{#1}}

\renewcommand{\O}{\mathcal{O}} 
 
\newcommand{\GP}{\mathcal{GP}}
\newcommand{\Id}{\vec{I}}

\newcommand{\II}{\mathbb{I}}

\newcommand{\y}{\vec{y}}
\newcommand{\x}{\vec{x}}

\usepackage{colonequals}

\usepackage{nicefrac}

\usetikzlibrary{arrows,shapes,plotmarks}

\tikzset{>=stealth'} 
\tikzstyle{graphnode} = 
   [circle,draw=black,minimum size=22pt,text centered,text
     width=22pt,inner sep=0pt] 
\tikzstyle{var}   =[graphnode,fill=white]
\tikzstyle{obs}   =[graphnode,fill=black,text=white]
\tikzstyle{fac}   =[rectangle,draw=black,fill=black!25,minimum size=5pt]
\tikzstyle{facprior} =[rectangle,draw=black,fill=black,text=white,minimum size=5pt]
\tikzstyle{edge}  =[draw=white,double=black,thick,-]
\tikzstyle{prior} =[rectangle, draw=black, fill=black, minimum size=
5pt, inner sep=0pt]
\tikzstyle{dirprior} = [circle, draw=black, fill=black, minimum
size=5pt, inner sep=0pt]

\DeclareSymbolFont{stmry}{U}{stmry}{m}{n}
\DeclareMathSymbol\leftarrowtriangle\mathrel{stmry}{"5E}
\DeclareMathSymbol\rightarrowtriangle\mathrel{stmry}{"5F}

\renewcommand{\to}{\operatorname*{\rightarrowtriangle}}

\newcounter{PHcomment}

\newcounter{MGcomment}

\newcounter{MOcomment}

\newlength\figheight
\newlength\figwidth

\begin{document}

\title{Probabilistic Numerics and Uncertainty in Computations}

\author{
Philipp Hennig$^{1}$, Michael A Osborne$^{2}$ and Mark Girolami$^{3}$}

\address{$^{1}$Max Planck Institute for Intelligent Systems,
    T\"ubingen, Germany\\
$^{2}$University of Oxford, United Kingdom\\
$^{3}$University of Warwick, United Kingdom}

\subject{Statistics, Computational Mathematics, Artificial Intelligence}


\corres{Philipp Hennig\\
\email{phennig@tue.mpg.de}}

\begin{abstract}  
We deliver a call to arms for \emph{probabilistic numerical methods}: algorithms for numerical tasks, including linear algebra, integration, optimization and solving differential equations, that return uncertainties in their calculations. Such uncertainties, arising from the loss of precision induced by numerical calculation with limited time or hardware, are important for much contemporary science and industry. Within applications such as climate science and astrophysics, the need to make decisions on the basis of computations with large and complex data has led to a renewed focus on the management of numerical uncertainty. We describe how several seminal classic numerical methods can be interpreted naturally as probabilistic inference. We then show that the probabilistic view suggests new algorithms that can flexibly be adapted to suit application specifics, while delivering improved empirical performance. We provide concrete illustrations of the benefits of probabilistic numeric algorithms on real scientific problems from astrometry and astronomical imaging, while highlighting open problems with these new algorithms. Finally, we describe how probabilistic numerical methods provide a coherent framework for identifying the uncertainty in calculations performed with a combination of numerical algorithms (e.g.~both numerical optimisers and differential equation solvers), potentially allowing the diagnosis (and control) of error sources in computations.  
\end{abstract}


\begin{fmtext}
\end{fmtext}

\maketitle

\section{Introduction}
\label{sec:introduction}

Probability theory is the quantitative framework for scientific inference \cite{jeffreys1961theory}. It codifies how observations (data) combine with modelling assumptions (prior distributions and likelihood functions) to give evidence for or against a hypothesis and values of unknown quantities. There is continuing debate about how prior assumptions can be chosen and validated \cite[e.g.][\textsection 1.2.3]{hutter10:_univer_artif_intel}; but the role of probability as the language of uncertainty is rarely questioned. That is, as long as the subject of inference is a physical variable. What if the quantity in question is a mathematical statement, the solution to a computational task? Does it make sense to assign a probability measure $p(x)$ over the solution of a linear system of equations $Ax=b$ if $A$ and $b$ are known? If so, what is the meaning of $p(x)$, and can it be identified with a notion of `uncertainty'?  If one sees the use of probability in statistics as a way to remove ``noise'' from ``signal'', it seems misguided to apply it to a deterministic mathematical problem. But noise and stochasticity are themselves difficult to define precisely. Probability theory does not rest on the notion of randomness (\emph{aleatory} uncertainty), but extends to quantifying \emph{epistemic} uncertainty, arising solely from missing information\footnote{Many resources discussing epistemic uncertainty can be found, at the time of writing, at \url{http://understandinguncertainty.org} (the authors are not affiliated with this web page)}. Connections between deterministic computations and probabilities have a long history. Erd\"os and Kac \cite{kac1940gaussian} showed that the number of distinct prime factors in an integer follows a normal distribution. Their statement is precise, and useful for the analysis of factorization algorithms \cite{knuth1976analysis}, even though it is difficult to ``sample'' from the integers. It is meaningful without appealing to the concept of epistemic uncertainty. Probabilistic and deterministic methods for inference on physical quantities have shared dualities from very early on: Legendre introduced the method of least-squares in 1805 as a deterministic \emph{best fit} for data without a probabilistic interpretation. Gauss' 1809 probabilistic formulation of the exact same method added a generative stochastic model for how the data might be assumed to have arisen. Legendre's least-squares is a useful method without the generative interpretation, but the Gaussian formulation adds the important notion of uncertainty (also interpretable as model capacity) that would later become crucial in areas like the study of dynamical systems.\footnote{In fact, the very same connection between least-squares estimation and Gaussian inference has been re-discovered repeatedly, simply because least-squares estimation has been re-discovered repeatedly after Legendre, under names like ridge regression \cite{ridge_regression}, Kriging \cite{krige1951statistical}, Tikhonov's method \cite{tikhonov1943stability}, and so on. The fundamental connection is that the normal distribution is the exponential of the square $\ell_2$ norm. Because the exponential is a monotonic function, minimizing an $\ell_2$-regularized $\ell_2$ loss is equivalent to maximizing the product of Gaussian prior and likelihood. In this sense, this paper is adding numerical mathematicians to the list as yet another group of re-discoverers of Gaussian inference. Ironically, this list includes Gauss himself, because Gaussian elimination, introduced in the very same paper as the Gaussian distribution itself~\cite{gauss1809theoria}, can be interpreted as a conjugate-direction method~\cite{hestenes1952methods}, and thus as Gaussian regression~\cite{2014arXiv1402.2058H}. See also~\cite{shental2008gaussian}.}

Several authors \cite{diaconis88:_bayes, ohagan92:_some_bayes_numer_analy,skilling1991bayesian} have shown that the language of probabilistic inference can be applied to numerical problems, using a notion of uncertainty about the result of an intractable or incomplete computation, and giving rise to methods we will here call \emph{probabilistic numerics}\footnote{The probabilistic numerics community web site can be found at \url{http://www.probabilistic-numerics.org}. }. In such methods, uncertainty regularly arises solely from the lack of information inherent in the solution of an ``intractable'' problem: A quadrature method, for example, has access only to finitely many function values of the integrand; an exact answer would, in principle, require infinitely many such numbers. Algorithms for problems like integration and optimization proceed iteratively, each iteration providing information improving a running estimate for the correct answer. Probabilistic numerics provides methods that, in place of such estimates, update probability measures over the space of possible solutions. As Diaconis noted \cite{diaconis88:_bayes}, it appears that Poincar\'e proposed such an approach already in the 19th century. The recent explosive growth of automated inference, and the increasing importance of numerics for science, has given this idea new urgency.

This article connects recent results, promising applications, and central questions for probabilistic numerics. We collate results showing that a number of basic, popular numerical methods can be identified with families of probabilistic inference procedures.  The probability measures arising from this new interpretation of established methods can offer improved performance, enticing new functionality, and conceptual clarity; we demonstrate this with examples drawn from astrometry and computational photography. The article closes by pointing out the propagation of uncertainty through computational pipelines as a guiding goal for probabilistic numerics.

It may be helpful to separate the issues discussed here clearly from other areas of overlap between statistics and numerical mathematics: We here focus on well-posed deterministic problems, identifying degrees of uncertainty arising from the computation itself. This is in contrast to the notion of uncertainty quantification \cite[e.g.][]{kennedy2001bayesian,kiureghian2009aleatory}, which identifies degrees of freedom in ill-posed problems, and where epistemic uncertainty arises from the set-up of the computation, rather than the computation itself. It also differs from several concepts of stochastic numerical methods, which use (aleatory) random numbers either to quantify uncertainty from repeated computations \cite[e.g.][]{iman1988investigation}, or to reduce computational cost through randomly chosen projections \cite[see for example][]{liberty2007randomized,halko2011finding}.

It is also important to note that, as a matter of course, existing frameworks already analyse and estimate the error created by a numerical algorithm. Theoretical \emph{analysis} of computational errors generally yields convergence rates---bounds up to an unknown constant---made under certain structural assumptions. Where probabilistic numerical methods are derived from classic ones, as described below, they naturally inherit such analytical properties. At runtime, the error is also \emph{estimated} for the specific problem instance. Such run-time error estimation is frequently performed by monitoring the dynamics of the algorithm's main estimate \cite[see][\textsection II.4, for ODEs]{hairer87:_solvin_ordin_differ_equat_i},~\cite[\textsection 4.5, for quadrature]{davis2007methods}. Similarly, in optimization problems, the magnitude of the gradient is often used to monitor the algorithm's progress. Such error estimates are informal, as are the solution estimates computed by the numerical method itself. They are meant to be used locally, mostly a criteria for the termination of a method, and the adaptation of its internal parameters. They can not typically be interpreted as a property (e.g.~the variance) of a posterior probability measure, and thus can not be communicated to other algorithms, and thus can not be embedded in a larger framework of error estimation. They also do not usually inform the design of the numerical algorithm itself; instead, they are a diagnostic tool added post-hoc. Below, we argue that the estimation of errors should be given a formal framework, and that probability theory is uniquely suited for this task. Describing numerical computations as inference on a latent quantity yields a joint, consistent, framework for the construction of solution and error estimates. The inference perspective can provide a natural intuition that may suggest extensions and improvements. And the probabilistic framework provides a \emph{lingua franca} for numerical computations, which allows the communication of uncertainty between methods in a chain of computation.

\section{Probabilistic numerical methods from classical ones}
\label{sec:design-prob-numer}

Numerical algorithms estimate quantities not directly computable, using the results of more readily available computations. Even existing numerical methods can thus be seen as inference rules, reasoning about latent quantities from ``observables'', or ``data''. At face value, this connection between inference and computation is vague. But several recent results have shown that it can be made rigorous, such that established deterministic rules for various numerical problems can be understood as maximum a-posteriori estimates under specific priors (hypothesis classes) and likelihoods (observation models). The recurring picture is that there is a one-to-many relationship between a classic numerical method for a specific task and a family of probabilistic priors which give the same maximum posterior estimate but differing measures of uncertainty. Choosing one member of this family amounts to fitting an uncertainty, a task we call \emph{uncertainty calibration}. The result of this process is a numerical method that returns a point estimate surrounded by a probability measure of uncertainty, such that the point estimate inherits the proven theoretical properties of the classic method, and the uncertainty offers new functionality.

\subsection{Quadrature}
\label{sec:quadrature}

\begin{figure}
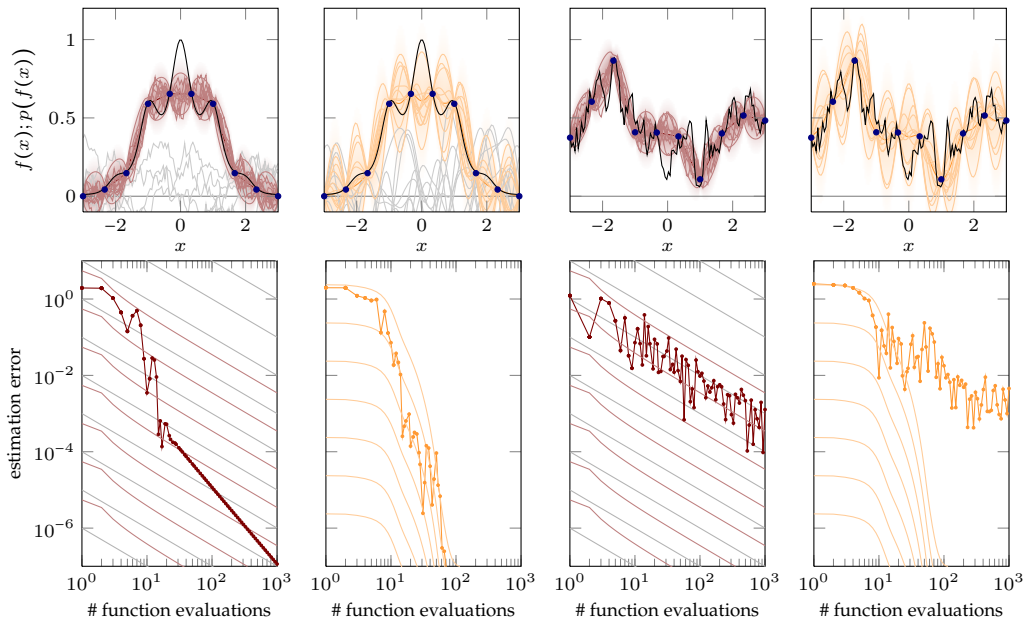

  \centering
  \setlength{\figwidth}{.2\textwidth}
  \setlength{\figheight}{0.2\textwidth}
  \scriptsize
  \mbox{
    \input{fig/quadrature_Wiener.tikz}\hspace{4mm}%
    \input{fig/quadrature_SE.tikz}\hspace{4mm}
    \input{fig/quadrature_Wiener_rough.tikz}\hspace{4mm}%
    \input{fig/quadrature_SE_rough.tikz}\hspace{4mm}%
    }\newline
  \mbox{\setlength{\figheight}{0.3\textwidth}%
    \input{fig/quadrature_convergence_W1.tikz}
    \input{fig/quadrature_convergence_SE.tikz}
    \input{fig/quadrature_convergence_W1_rough.tikz}
    \input{fig/quadrature_convergence_SE_rough.tikz}
    }
  \caption{\label{fig:quad} Quadrature rules, illustrating the challenge of uncertainty-calibration. Left two columns: Top row: Function $f(x)$ (black line), is approximately integrated using two different Gaussian process priors (left: linear spline; right: exponentiated quadratic), giving posterior distributions and mean estimates. Gray lines are functions sampled from the prior. The thick coloured line is the posterior mean, thin lines are posterior samples, and the delineation of two marginal standard deviations. The shading represents posterior probability density. Bottom row: As the number of evaluation points increases, the posterior mean (thick line with points) converges to the true integral value; note the more rapid convergence of the exponentiated quadratic prior. The posterior covariance provides an error estimate whose scale is defined by the posterior mean alone (each thin coloured line in the plots corresponds to a different instance of such an estimate). But it is only a meaningful error estimate if it is matched well to the function's actual properties. The left plot shows systematic difference between the convergence of the real error and the convergence of the estimated error under the linear spline, whereas convergence of the estimated error under the exponentiated quadratic prior is better calibrated to the real error. Grey grid lines in the background, bottom left, correspond to $\O(N^{-1})$ convergence of the error in the number $N$ of function evaluations. Right two columns: The same experiment repeated with a function $f$ drawn from the spline kernel prior. For this function, the trapezoid rule is the optimal statistical estimator of the integral (note well-calibrated error measure in bottom row, third panel), while the Gaussian kernel GP is strongly over-confident.} \end{figure}

We term the Probabilistic Numeric approach to quadrature \emph{Bayesian quadrature.} Diaconis \cite{diaconis88:_bayes} may have been first to point out a clear connection between a Gaussian process regression model and a deterministic quadrature rule, an observation subsequently generalised by Wahba \cite[][\textsection 8]{wahba1990spline} and O'Hagan \cite{o1991bayes}, and also noted by \citet{BZMonteCarlo}. Details can be found in these works; here we construct an intuitive example highlighting the practical challenges of assigning uncertainty to the result of a computation. For concreteness, consider $f(x) = \exp\left[ -\sin^2(3x) - x^2 \right]$ (black in Figure \ref{fig:quad}, top). Evidently, $f$ has a compact symbolic form and $f(x)$ can be computed for virtually any $x\in\Re$ in nanoseconds. It is a wholly deterministic object. Nevertheless, the real number
\begin{equation}
  \label{eq:2}
  F = \int_{-3} ^3 f(x)\, \mathrm{d}x
\end{equation}
has no simple analytic value, in the sense that it can not be natively evaluated in low-level code. Quadrature rules offer ``black box'' estimates of $F$. These rules have been optimized so heavily \cite[e.g.][]{davis2007methods} that they could almost be called ``low-level'', but their results do not come with the strict error bounds of floating-point operations; instead, assumptions about $f$ are necessary to bound error. Perhaps the simplest quadrature rule is the trapezoid rule, which amounts to linear interpolation of $f$ (red line in Figure \ref{fig:quad}, top left): Evaluate $f(x_i)$ on a grid $-3= x_1<x_2<\dots<x_N=3$ of $N$ points, and compute
\begin{equation}
  \label{eq:3}
  \hat{F}_{\text{midpoint}} = \sum_{i=2}
  ^N \frac{1}{2}\left[f(x_i)+f(x_{i-1})\right] (x_i-x_{i-1}).
\end{equation}

\subsection{Bayes-Hermite Quadrature} \label{sec:bayes-herm-quadr} A probabilistic description of $F$ requires a joint probability distribution over $f$ and $F$. Since $f$ lies in an infinite-dimensional Banach or Hilbert space, no Lebesgue measure can be defined. But Gaussian process measures can be well-defined over such spaces, and offer a powerful framework for quadrature \cite{wahba1990spline,RasmussenWilliams}. In particular, we may choose to model $f$ on $[-3,3]$ by $p(f)=\GP(f;0,k)$, a Gaussian process with vanishing mean $\mu=0$ and the linear spline covariance function \cite{minka2000deriving} (a stationary variant of the Wiener process)
\begin{equation}\label{eq:spline}
 k(x,x')=c(1+b-\nicefrac{1}{3}\,b\,|x-x'|), \qq \text{ for some }c,b>0. 
\end{equation}
More precisely, this assigns a probability measure over the function values $f(x)$ on $[-3,3]$, such that any restriction to finitely many evaluations $\y=[f(x_1),\dots,f(x_M)]$ at $X=[x_1,\dots,x_N]$ is jointly Gaussian distributed with zero mean and covariance $\cov \bigl[f(x_i),f(x_j)\bigr] = k(x_i,x_j)$. Samples drawn from this process are shown in grey in the top left of Figure \ref{fig:quad}. These sample paths represent the hypothesis class associated with this Gaussian process prior: They are continuous, but not differentiable, in mean square expectation \cite[][\textsection~2.2]{adler1981geometry}. Because Gaussian processes are closed under linear projections \cite[e.g.][p.~27]{adler1981geometry}, this distribution over $f$ is identified with a corresponding univariate Gaussian distribution \cite{minka2000deriving}, 
\begin{equation}
  \label{eq:4}
  p(F) = \N\left[ F;0,\iint_{-3} ^3
  k(\tilde{x},\tilde{x}')\,\mathrm{d}\tilde{x}\,\mathrm{d}\tilde{x}' 
    = c(1 + \nicefrac{b}{3}) \right],
\end{equation}
on $F$. The strength of this formulation is that it provides a clear framework under which observations $f(x_i)$ can be incorporated. The measure on $p(F)$ conditioned on the collected function values is another scalar Gaussian distribution
 \begin{multline}
   \label{eq:5}
     p \bigl(F\g \y \bigr) = \N\left[F;\int_{-3} ^3 k(\tilde{x},X)K^{-1}\y\,\mathrm{d}\tilde{x},
       \iint_{-3} ^3 k(\tilde{x},\tilde{x}') - k(\tilde{x},X)K^{-1}k(X,\tilde{x}')\,
        \mathrm{d}\tilde{x}\,\mathrm{d}\tilde{x}'\right],
 \end{multline}
 where $k(\tilde{x},X)=[k(\tilde{x},x_1),\dots,k(\tilde{x},x_N)]$ and $K$ is the symmetric positive definite matrix with elements $K_{ij} = k(x_i,x_j)$. The final expression in the brackets, the posterior covariance can be optimized with respect to $X$ to \emph{design} a ``most informative'' dataset. For the spline kernel~\eqref{eq:spline}, this leads to placing $X$ on a regular, equidistant, grid \cite{minka2000deriving}. 

 The spline kernel $k$ of \eqref{eq:spline} is a piecewise linear function with a single point of non-differentiability at $x=x'$. The posterior mean $k(\tilde{x},X)K^{-1}\y$ is a weighted sum of $N$ such kernels centred on the evaluation points $X$ and constrained by the likelihood to pass through the values $\y$ at the nodes $X$. Thus, the posterior mean is a linear spline interpolant of the evaluations, and the posterior mean of Equation~\eqref{eq:5} is exactly equal to the trapezoid-rule estimate from $X$ \cite[see also][]{diaconis88:_bayes}. That is, Bayesian quadrature with a linear spline prior provides a probabilistic interpretation of the trapezoidal rule; it supplements the estimate with a full probability distribution characterising uncertainty. The corresponding conditional on $f$ is a Gaussian process whose mean is the piecewise linear red function in the top left of Figure \ref{fig:quad} (the Figure also represents the conditional distribution $p \bigl(f\g f(x_1),\dots,f(x_N)\bigr)$ as a red cloud, and some samples). Various other, more elaborate quadrature rules (e.g.\ higher order spline interpolation and Chebyshev polynomials) can be cast probabilistically in analogous ways, simply by changing the covariance function $k$ \cite{wahba1990spline,o1991bayes,diaconis88:_bayes}.

 \subsection{Added value, and challenge, of probabilistic output}
\label{sec:added-value-chall}

A non-probabilistic analysis of the trapezoid-rule is that, for sufficiently regular functions $f$, this rule converges at a rate of at least $\O(N^{-1})$ to the correct value for $F$ as $N$ increases \cite[e.g.][Eq.~2.1.6]{davis2007methods} (other quadrature rules, like cubic splines, have better convergence for differentiable functions). Figure \ref{fig:quad} bottom left shows the mean estimate (equivalently, the trapezoid-rule) as a thick line. Grey help lines represent $\O(N^{-1})$ convergence. Clearly, the asymptotic behaviour is steeper than those lines.\footnote{This, too, is a well-known result: If the integrand is differentiable, rather than just continuous, the trapezoid rule has quadratic convergence rate \cite[][Eq.~2.1.12]{davis2007methods}.} So the theoretical bound is correct, but it is also of little practical value: No linear bound is tight from start to convergence.

On the probabilistic side, the standard deviation of the conditional (\ref{eq:5}) is regularly interpreted as an estimate of the imprecision of the mean estimate. In fact, this error estimate (thin red lines in Figure \ref{fig:quad}, top left) is analogous to the deterministic linear convergence bound for continuous functions. There is a family of such error bounds associated with the same posterior mean, with each line corresponding to a different value for the unknown constant in the bound and different values for the parameters $b$ and $c$ in the covariance.\footnote{ The initial behaviour of the red lines in the figure is a function of the scale $b$, which relates to assumptions about the rate at which the asymptotic quadratic convergence is approached.} It may seem as though the probabilistic interpretation had added nothing new, but since this view identifies quadrature rules of varying assumptions as parameter choices in Gaussian process regression, it embeds seemingly separate rules in a hierarchical space of models, from which models with good error modelling can be selected by hierarchical probabilistic inference. This can be done without collecting additional data-points \cite[{\textsection}5]{RasmussenWilliams}.

More generally, the probabilistic numeric viewpoint provides a principled way to manage the parameters of numerical procedures. Where Markov Chain Monte Carlo procedures might require the hand-tuning of parameters such as step sizes and annealing schedules, Bayesian quadrature allows the machinery of statistical inference procedure to be brought to bear upon such parameters. A practical example of the benefits of approximate Bayesian inference for the (hyper-)parameters of a Bayesian quadrature procedure is given by \cite{osborne2012active}.

The function $f$ used in this example is much smoother than typical functions under the Gaussian process prior distribution associated with the trapezoid rule (shown as thin grey samples in the top left plot). The second column (using orange) of Figure~\ref{fig:quad} shows analogous experiments with the exponentiated quadratic covariance function $k(x,x') = \theta^2 \exp(-(x-x')^2/\lambda^2)$, corresponding to a very strong smoothness assumption on $f$ \cite{van2011information} (see grey samples in top right plot), giving a very quickly converging estimate. In this case, the ``error bars'' provided by the standard deviation converge in a qualitatively comparable way. 

Of course, the faster convergence of this quadrature rule based on the exponentiated quadratic-covariance prior is not a universal property. It is the effect of a much stronger set of prior assumptions. If the true integrand is rougher than expected under this prior, the quadrature estimate arising from this prior can be quite wrong. The right two columns of Figure~\ref{fig:quad} show analogous experiments on an integrand  that is a true sample from a Gaussian process with the spline covariance~\eqref{eq:spline}. In this case, the spline prior is the optimal statistical estimator by construction, and its error estimate is perfectly calibrated (third column of Figure~\ref{fig:quad}), while the exponentiated quadratic-kernel gives over-confident, and inefficient estimates (right-most column). 

Identifying the optimal regression model from a larger class, just based on the collected function values, requires more computational work than to fix a regressor from the start. But it also gives better calibrated uncertainty. Contemporary general-purpose quadrature implementations \cite[e.g.][]{davis2007methods} remain lightweight by recursively re-using previous computations. The above experiments show that it is possible to design Bayesian quadrature rules with well-calibrated posterior error estimates, but it remains a question how small the computational overhead from a probabilistic computation over these methods can be made. Even so, formulating quadrature as probabilistic regression precisely captures a trade-off between \emph{prior assumptions} inherent in a computation and the computational effort required in that computation to achieve a certain precision. Computational rules arising from a strongly constrained hypothesis class can perform much better than less restrictive rules \emph{if the prior assumptions are valid}. In the numerical setting---in contrast to many empirical situations in statistics---it is often possible to precisely check whether a particular prior assumption is valid: The machine performing the computation has access, at least in principle, to a formal, complete, description of its task, in form of the source code describing the task (the integrand, the optimization objective, etc.). Using this source code, it is for example to test, at runtime, whether, and how many times, an integrand is continuously differentiable~\cite[e.g.][]{Griewank2000EDP}. Using this information will in future allow the design of improved quadrature methods. Once one knows that general purpose quadrature methods effectively use a Gaussian process prior over function values, it is natural to ask whether this prior actually incorporates the salient information available for one's specific problem. Including such information in the prior leads to customized, ``tailored'', numerical methods that can perform better.

Probabilistic inference also furnishes the framework required to tackle numerics tasks using decision theory. For quadrature, the decision problem to be solved is that of node selection: that is, the determination of the optimal positions for points (or nodes) at which to evaluate the integrand. With the definition of an appropriate loss function, such as the posterior variance of the integral, such nodes can be optimally selected by minimising the expected loss function. It seems clear that this approach can improve upon the simple uniform grids of many traditional quadrature methods, and can enable active learning: where surprising evaluations can influence the selection of future nodes.

This decision-theoretic approach stands in contrast with that adopted by randomised, Monte Carlo, approaches to quadrature. Such approaches generate nodes with the use of a pseudo-random number generator, injecting additional epistemic uncertainty (about the value of the generator's outputs) into a procedure designed to reduce the uncertainty in an integral. It is worth noting that the use of pseudo-random generators burdens the procedure with additional computational overhead: pseudo-random numbers are cheap, but not free. The principal feature of Monte Carlo approaches is their conservative nature: the Monte Carlo policy will always, eventually, take an additional node arbitrarily close to an existing node. The disadvantage of this strategy is its waste of valuable evaluations: the convergence rate of Monte Carlo techniques, $\O(N^{-1/2})$, is clearly improved upon by both traditional quadrature and Bayesian quadrature methods. This problem is only worsened by the common discarding of evaluations known as `burn-in' and `thinning'. The advantage of Monte Carlo, of course, is its robustness to even highly non-smooth integrands. However, Bayesian quadrature can realise more value from evaluations by exploiting known structure (e.g. smoothness) in the integrand.

\subsection{Empirical evaluation of Bayesian quadrature for astrometry} 
\label{sub:empirical_evaluation_of_bayesian_quadrature_for_astrometry}

We illustrate this on an integration problem drawn from astrometry, the measurement of the motion of stars.  In order to validate astrometric analysis, we aim to recover the number of planets present in synthetic data generated (with a known number of planets) to mimic that produced by an astrometric facility such as the GAIA satellite \cite{sozzetti2014astrometric}.\footnote{The authors are grateful to H. Parviainen and S. Aigrain for providing data and code examples.}  Here, quadrature's task is to compute the model evidence (marginal likelihood) $Z=\int p(\data \mid \theta, \model)~p(\theta \mid \model )\,\mathrm{d}\theta$ of a model $\model$ for orbital motions, where $\data$ are the gathered observations and $\theta$ are the model parameters. Specifically, it is of interest to compare the evidence for models including differing numbers of exoplanets. 
For the following example, to provide a focal problem upon which to compare quadrature methods, we compute the evidence of a model with two such planets on data generated with a two-planet model. The corresponding integral is analytically intractable, with a multi-modal integrand (the likelihood $p(\data \mid \theta, \model)$) and a 19-dimensional $\theta$ rendering the quadrature problem challenging.

\begin{figure}
  \centering\scriptsize
  \setlength{\figwidth}{.9\textwidth}
  \setlength{\figheight}{0.15\textwidth}
  \mbox{\input{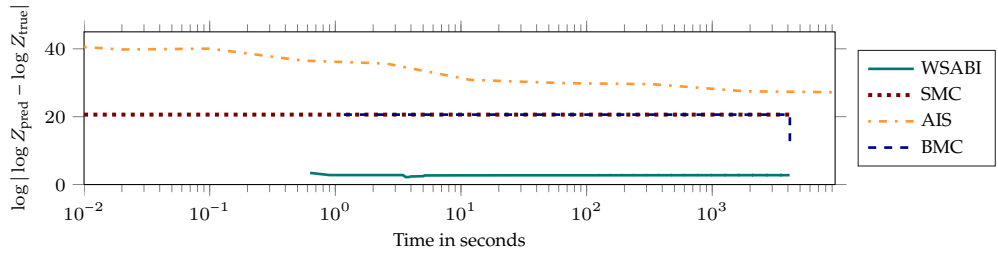}}
  \caption{
    We compare different quadrature methods in computing the 19-dimensional integral giving model evidence in an astrometry application. 
    The Bayesian quadrature algorithm (WSABI) that employs active selection of nodes, along with prior knowledge of the smoothness and non-negativity of the integrand, converges faster than the Monte Carlo approaches: simple Monte Carlo (SMC) and Annealed Importance Sampling (AIS). 
    Note that Bayesian Monte Carlo (BMC) is a Bayesian quadrature algorithm that uses the same samples as SMC, explaining their similar performance.
    Note also the performance improvement offered by WSABI over BMC, suggesting the crucial role played by the active selection of nodes.
    }
\label{fig:WSABI}
\end{figure}

As the probabilistic method, we use a recent algorithm: Warped, Sequential, Active Bayesian Integration (WSABI\footnote{Matlab code for WSABI is available at {\tt https://github.com/OxfordML/wsabi}.}) \cite{gunter2014sampling}. 
WSABI is a Bayesian quadrature algorithm that uses an internal probabilistic model that is well-calibrated to the suspected properties of the problem's integrand. 
Firstly, like the exponentiated quadratic, WSABI's covariance function is suitable for smooth integrands, as are expected for the problem.
Secondly, and beyond what is achievable with classic quadrature rules, WSABI explicitly encodes the fact that the integrand (the likelihood $p(\data \mid \theta, \model)$) is strictly positive.
WSABI also makes use of a further opportunity afforded by the probabilistic numeric approach: it actively select nodes so as to minimise the uncertainty in the integral. 
This final contribution permits nodes to be selected that are far more informative than gridded or randomly selected evaluations.
We compare this algorithm against two different Monte Carlo approaches to the problem: Annealed Importance Sampling (AIS) \cite{neal2001annealed} (which was implemented with a Metropolis--Hastings sampler) and simple Monte Carlo (SMC). 
We additionally compared against Bayesian Monte Carlo (BMC) \cite{BZMonteCarlo}, a Bayesian quadrature algorithm using the simpler exponentiated quadratic model, and whose nodes were taken from the same samples selected by SMC. 
``Ground truth'' ($Z_{\text{true}}$) was obtained through exhaustive SMC sampling ($10^6$ samples). 
The results in Figure \ref{fig:WSABI} show that the probabilistic quadrature method achieves improved precision drastically faster than the Monte Carlo estimates. 
It is important to point out that the plot's abscissa is ``wall-clock'' time, not algorithmic steps.  Probabilistic algorithms need not be expensive.

\subsection{Linear Algebra}
\label{sec:linear-algebra}

Computational linear algebra covers various operations on matrices. We will here focus on linear optimization, where $\vec{b}\in\Re^{N}$ is known and the task is to find $\x\in\Re^N$ such that $A\x = \vec{b}$ where $A\in\Re^{N\times N}$ is symmetric positive definite.  We assume access to projections $A\vec{s}$ for arbitrary $\vec{s}\in\Re^N$. If $N$ is too large for exact inversion of $A$, a widely known approach is the method of conjugate gradients (CG) \cite{hestenes1952methods}, which produces a convergent sequence $\x_0,\dots,\x_M$ of improving estimates for $\x$. Each iteration involves one matrix-vector multiplication and a small number of linear operations, to produce the update $\vec{s}_i=\x_{i}-\x_{i-1}$, and an ``observation'' $\y_i=A\vec{s}_i$. CG's good performance has been analysed extensively \cite[e.g.][\textsection 5.1]{nocedal1999numerical}. 

We will use the shorthands $S_M=[\vec{s}_1,\dots,\vec{x}_M]$ and $Y_M=[\y_1,\dots,\y_M]$ for the set of projection directions and ``observations'' after $M$ iterations of the method. Defining a probabilistic numerical algorithm requires a joint probability measure $p(Y_M,A,\vec{b},\x\g S_M)$ over all involved variables, conditioned on the algorithm's active design decisions; and an ``action rule'' (design, policy, control law) describing how the algorithm should collect data. Recent work by Hennig \cite{2014arXiv1402.2058H} describes such a model which exactly reproduces the sequence $\{\x_i\}_{i=1,\dots,M\leq N}$ of CG: As prior, choose a Gaussian measure over the (assumed to exist) inverse $H=A^{-1}$ of $A$ 
\begin{equation}
  \label{eq:7}
  p(H) = \N(\vect{H};\vect{H}_0,\Gamma (W\otimes W)\Gamma\Trans).
\end{equation}
This implies a Gaussian prior over $x=Hb$, and a non-Gaussian prior over $A$. In (\ref{eq:7}), $\N$ is a Gaussian distribution over the $N^2$ elements of $H$ stacked into a vector $\vect{H}$. $\Gamma$ is the symmetrization operator $(\Gamma \vect{A} = \nicefrac{1}{2}(\vect{A+A \Trans}))$, and $\otimes$ is the Kronecker product ($\Gamma(W\otimes W)\Gamma\Trans$ is also known as the symmetric Kronecker product \cite{loan00:_kronec}). As projections are presumed noise-free, the likelihood is the Dirac distribution, the limit of a Gaussian with vanishing width, which can also be seen as performing a strict conditioning of the prior on the observed function values,
\begin{equation}
  \label{eq:8}
  \begin{split}
    p(Y_M\g A,\vec{b},\x,S_M) &= \prod_i \delta(A\vec{s}_i - \y_i) = \prod_i \delta(\vec{s}_i -
 H\y_i).
  \end{split}
\end{equation}
The action rule at iteration $i$ is to move to $\x_{i+1} = \x_{i} - \alpha \hat{H}_{i} (A\x_{i} - \vec{b})$, where $\hat{H}_{i}$ is the posterior mean $p(H\g Y_M)$. The optimal step $\alpha$ can be computed exactly in a linear computation.

Eq.~(\ref{eq:7}) requires choices for $\vect{H}_0$ and $W$. The prior mean is set to unit, $\vect{H}_0=\Id$. If $W\in\Re^{N\times N}$ is positive definite, then $p(H)$ assigns finite measure to every symmetric $H\in\Re^{N\times N}$, and the algorithm converges to the true $\x$ in at most $N$ steps, assuming exact computations \cite{2014arXiv1402.2058H}.  In this sense it is non-restrictive, but of course some matrices are assigned more density than others. If $W$ is set to a value among the set $\{W\in\Re^{N\times N}\g W \text{ symm. positive definite and } WY_M=\sigma S_M, \sigma\in\Re_+\}$ (this includes the true matrix $W=H$, but does not require access to it), then the resulting algorithm \emph{exactly} reproduces the iteration sequence $\{ \x_i\}_{i=0,\dots,N}$ of the CG method, and can be implemented in the exact same way. However, this derivation also provides something new: A Gaussian posterior distribution $p(H\g \y)=\N(H;H_M,\Gamma(W_M\otimes W_M)\Gamma\Trans)$ over $H$. Details can be found in \cite{2014arXiv1402.2058H}, for the present discussion it suffices to know that both the posterior mean $H_M$ and covariance parameter $W_M$ are of manageable form. In particular, $H_M=\Id+ UEU\Trans$ with a diagonal matrix $E\in\Re^{2M\times 2M}$ and ``skinny'' matrices $U\in\Re^{N\times 2M}$, which are a function of the steps $\vec{s}$ and observed projections $\vec{y}$. So, as in the case of quadrature, there is a family of Gaussian priors of varying width (scaled by $\sigma$), such that all members of the family give the same posterior mean estimate. And this posterior mean estimate is identical to a classic numerical method (CG). But each member of the family gives a different posterior covariance---a different uncertainty estimate.

It is an interesting question to which degree the uncertainty parameter $W_M$ can be designed to give a meaningful error estimate. Some answers can be found in \cite{2014arXiv1402.2058H}. Interestingly, the equivalence-class of prior covariances $W_0$ that match CG in the mean estimate has more degrees of freedom than the number $M$ of observations (matrix-vector multiplications) collected by the algorithm during its typical runtime. Fitting the posterior uncertainty $W_M$ thus requires strong regularization. The method advocated in \cite{2014arXiv1402.2058H} constructs such a regularized estimator exclusively from scalar numbers already collected during the run of the CG method, thus keeping computational overhead very small.

But, as in quadrature, there are valuable applications for the probabilistic formulation that do not strictly require a well-calibrated \emph{width} of the posterior. Applications that make primary use of the posterior mean may just require \emph{algebraic structure} in the prior, up to an arbitrary scaling constant, to incorporate available, helpful, information. We will not do this here and instead highlight another use of probabilities over point estimates: the propagation of knowledge from one linear problem to another related problem. 

The approach described in the following is known in the numerical linear algebra community as the \emph{recycling of Krylov sequences} \cite{doi:10.1137/040607277}. 
However, while the framework of classic numerical analysis required challenging and bespoke derivation of this result, it follows naturally in the probabilistic viewpoint as the extension of a parametric regressor to a filter on a time-varying process. 
The probabilistic formulation of computation uses the universal and unique language of inference to enable the solution of similar problems across the breadth of numerics using similar techniques. This is in contrast to the compartmentalised state of current numerics, which demands distinct expert knowledge of each individual numeric problem in order to make progress towards it. 
Further, there is social value in making results accessible to any practitioner with a graduate knowledge of statistics.

\begin{figure*}
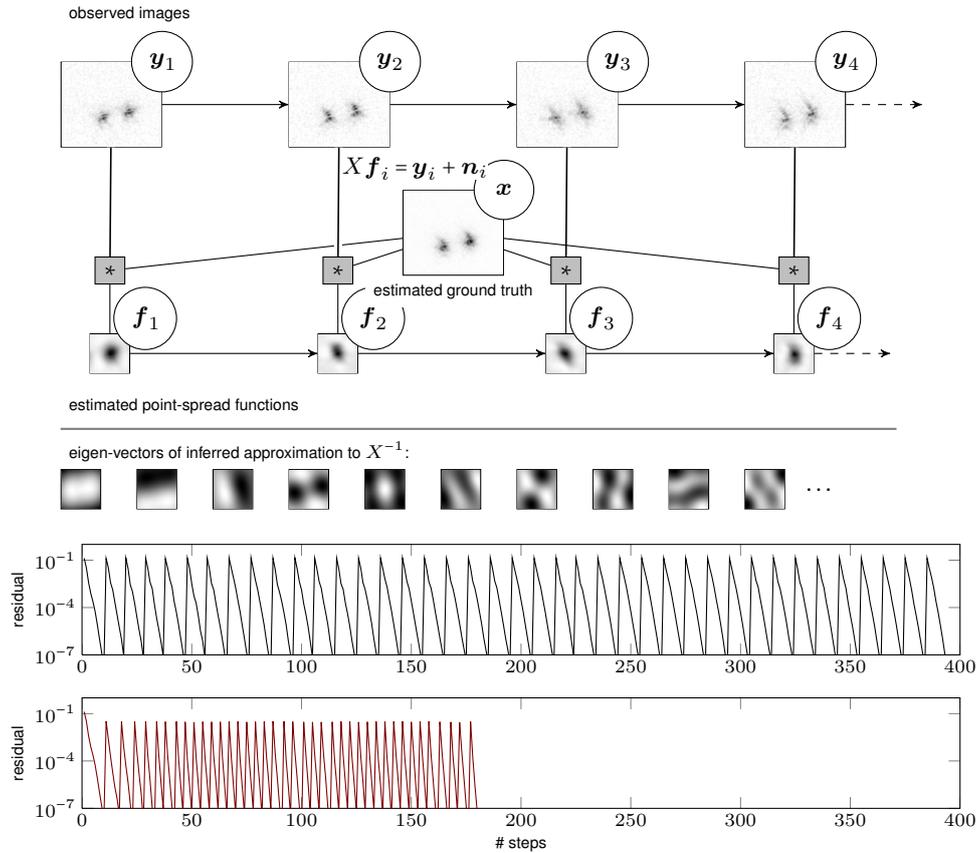

  \centering
  \mbox{\input{fig/LinearExp.tikz}}\\[\bigskipamount]
  \setlength{\figwidth}{.9\textwidth}%
  \setlength{\figheight}{3.5cm}%
  \scriptsize \sf
  \input{fig/batch_gradients_run}
  \caption{\label{fig:CG} Solving sequences of linear problems by utilizing the probabilistic interpretation of conjugate gradients, also known as the recycling of Krylov sequences. Top: A sequence of observed astronomical images $\y_i$ is modelled as the convolution of a stationary true image $\x$ with a time-varying point-spread function $\vec{f}_i$ (the matrix $X$ encodes the convolution with $\x$). Each individual deconvolution task requires one run of a linear solver, here chosen to be the method of conjugate gradients. Bottom plots: If each problem is solved independently, each instance of the solver progresses similarly (top plot. Optimization residuals over time. Each ``jump'' in the residual is the start of a new frame/deconvolution problem). If the posterior mean implied by the probabilistic interpretation of the solver is communicated from one problem to the next, the solvers progress increasingly faster (bottom plot, note decreasing number of steps for each deconvolution problem, and decrease, by about one order of magnitude, of initial residuals). Middle row: Over time, the vectors spanning this posterior mean for $X^{-1}$ converge to a generic basis for point-spread functions.}
\end{figure*}

In many set-ups, the same $A$ features in a sequence of problems $A\x_t=\vec{b}_t$, $t=1,\dots$. In others, a map $A_t$ changes slightly from step $t$ to step $t+1$. Figure \ref{fig:CG} describes a blind deconvolution problem from astronomical imaging:\footnote{The authors thank S. Harmeling, M. Hirsch, S. Sra and B. Sch\"olkopf \cite{harmeling2009online} for providing data.} The task is to remove an unknown linear blur from a sequence $\y_1,\dots,\y_K$ of astronomical images 
\cite{harmeling2009online}. Atmospherical disturbances create a blur that continuously varies over time. The model is that each frame $\y_i$ is the noisy result of convolution of the same ground truth image $\x$ with a spatially varying blur kernel $\vec{f}_i$, i.e. $\y_i = \vec{f}_i * \x + \vec{n}_i$, where $\vec{n}_i$ is white Gaussian noise. A matrix $X$ encodes the convolution operation as $\y_i=X\vec{f}_i+\vec{n}_i$ (it is possible to ensure $X$ is positive definite). A blind deconvolution algorithm iterates between estimating $X$ and estimating the $\vec{f}_i$ \cite{harmeling2009online}. Each iteration thus requires the solution of $K$ linear problems (with fixed $X$) to find the $\vec{f}_i$, then one larger linear problem to find $\x$. The na\"ive approach of running separate instances of CG wastes information, because the $K$ linear problems all share the matrix $X$, and from one iteration to the next, the matrix $X$ changes less and less as the iterations approach convergence. Instead, information can be \emph{propagated} between related computations using the probabilistic interpretation of CG, by starting each computation in the sequence with a prior mean $H_0$ set to the posterior mean $H_M$ of the preceding problem. Due to the low-rank structure of $H_M$, this has low cost. To prevent a continued rise in computational costs as more and more linear problems are solved, $H_M$ can be restricted to a fixed rank approximation after each inner loop, also at low cost. 
The plots in the lower half of Figure~\ref{fig:CG} show the increase in computational efficiency: For the baseline of solving 40 linear problems independently in sequence, each converges about equally fast (top plot, each ``jump'' is the start of a new problem). The lower plot shows optimization progress of the same 40 problems when information is propagated from one problem to the next. The first problem amounts to standard CG, while subsequent iterations can make increasingly better use of the available information. The Figure also shows the dominant eigenvectors of the inferred posterior means of $X^{-1}$ after $K=40$ subsequent linear problems, which converge to a relatively generic basis for point-spread functions. Although not strictly correct, this scheme can be intuitively understood as \emph{inferring a pre-conditioner} across the sequence of problems, by Bayesian filtering (the technical caveat is that the described scheme does not re-scale the linear problem itself, as a pre-conditioner would, it just shifts the initial solution estimate). 

\subsubsection{Further Areas}
\label{sec:remark}

  These examples highlight the areas of quadrature and linear algebra. Analogous results, identifying existing numerical methods with maximum a-posteriori estimators have been established in other areas, too. We do not experiment with them here, but they help complete the picture of numerical methods as inference across problem boundaries:

In non-linear optimization, quasi-Newton methods like the BFGS rule are deeply connected to conjugate gradients (in linear problems, BFGS and CG are the same algorithm \cite{nazareth1979relationship}). The BFGS algorithm can be interpreted as a specific kind of autoregressive generalization of the Gaussian model for conjugate gradients \cite{hennig13:_quasi_newton_method}. Among other things, this allows explicit modelling of noise on evaluated gradients, a pressing issue in large-scale machine learning \cite{StochasticNewton}.

One currently particularly exciting area for probabilistic numerics is ordinary differential equations, more specifically the solution of initial value problems (IVPs) of the form $\nicefrac{d\x}{dt}(t)=f(\x(t),t)$, where $x:\Re_+\to\Re^N$ is a real-valued curve parametrised by $t$ known to start at the initial value $\x(t_0)=\x_0$. Explicit Runge-Kutta methods are a basic and well-studied tool for such problems \cite{hairer87:_solvin_ordin_differ_equat_i}---while not necessarily state of the art, they nevertheless perform well on many problems, and are conceptually very clear. From the inference perspective, Runge-Kutta methods are linear extrapolation rules. At increasing nodes $t_0<t_{1}<\dots<t_i<t_s$ they repeatedly construct ``estimates'' $\hat{\x}(t_i)\approx \x(t_i)$ for the true solution which is used to collect an ``observation'' $\y_i=f\bigl(\hat{\x}(t_i),t_i \bigr)$, such that the estimate is a linear combination of previous observations:
\begin{equation} 
    \label{eq:9} 
    \hat{\x}(t_i) = \x_0 + \sum_{j<i} w_{ij}\y_j. 
  \end{equation} 
  Crucially, the weights $w_{ij}$ are chosen such that, after $s$ evaluations, the estimate has high convergence order: $\|\x(t_s) - \hat{\x}(t_s)\| = \mathcal{O}(h^p)$, with an order $p\leq s$ (typically, $p\leq 5$).

It is important to point out the central role that linear extrapolation, and linear computations more generally, play again here, as they did in the other numerical settings discussed above. It is not an oversimplification to note that numerical methods often amount to efficiently projecting an intractable problem into a tractable linear computation. Since the Gaussian family is closed under all linear operations, it is, perhaps, no surprise that Gaussian distributions play a central role in probabilistic re-interpretations of existing numerical methods. In the case of initial value problems, Gaussian process extrapolation was previously suggested by Skilling \cite{skilling1991bayesian} as a tool for their solution. Starting from scratch, Skilling arrived at a method that shares the linear structure of Eq. (\ref{eq:9}), but does not have the strong theoretical underpinning of Runge-Kutta methods, in particular, Skilling's method does not share the high convergence order of Runge-Kutta methods. But the probabilistic formulation allows for novel theoretical analysis of its own \cite{o.13:_bayes_uncer_quant_differ_equat}, and new kinds of applications, such as the marginalization over posterior uncertainty in subsequent computational steps and finite prior uncertainty over the initial value \cite{HennigAISTATS2014}. A considerably vaguer but much earlier observation was made, on the side of numerical mathematics, by Nordsieck \cite{nordsieck1962numerical}, who noted that a class of methods he proposed, and which were subsequently captured in a wider nomenclature of ODE solvers, bore a resemblance to linear electrical filters. These, in turn, are closely connected to Gaussian process regression through the notion of Markov processes.

Recently, Schober et al.~\cite{runge-kutta14:_schob_m} showed that these connections between Gaussian regression and the solution of IVP's, hitherto only conceptual, can be made tight. There is a family of Gauss-Markov priors that, used as extrapolation rules, give posterior Gaussian processes whose mean function exactly matches members of the Runge-Kutta family. (Thanks to its Markov property, the corresponding inference method can be implemented as a signal filter, and thus in linear computational complexity, like Runge-Kutta methods). Hence, as in the other areas of numerics, there is now a family of methods that returns the trusted point estimates of an established method, while giving a new posterior uncertainty estimate allowing new functionality.

\section{Discussion}
\label{sec:discussion}

\subsection{General Recipe for Probabilistic Numerical Algorithms}
\label{sec:gener-recipe-prob}

These recent results, identifying probabilistic formulations for classic numerical methods, highlight a general structure. Consider the problem of approximating the intractable variable $z$, if the algorithm has the ability to choose `inputs' $\x=\{x_i\}_{i=1,\dots}$ for computations that result in numbers $\y(\x)=\{y_i(x_i)\}_{i=1,\dots}$. A blueprint for the definition of probabilistic numerical methods requires two main ingredients:

\begin{enumerate}
\item A \emph{generative model} $p(z,\y(\x))$ for all variables involved---that is, a joint probability measure over the intractable quantity to be computed, and the tractable numerical quantities computed in the process of the algorithm. Like all (sufficiently structured) probability measures, this joint measure can be written as
\begin{equation}
    \label{eq:1}
    p\bigl(z,\y(\x)\bigr) = p(z)\,p\bigl(\y(\x)\big\g z\bigr),
  \end{equation}
  i.e. separated into a \emph{prior} $p(z)$ and a \emph{likelihood} $p\bigl(\y(\x)\big\g z\bigr)$. The prior encodes a hypothesis class over solutions, and assigns a typically non-uniform measure over this class. The likelihood explains how the collected tractable numbers $\y$ relate to $z$. It has the basic role of describing the numerical task. Often, in classic numerical problems, the likelihood is a deterministic conditioning rule, a point measure.
\item A \emph{design}, \emph{action} rule, or \emph{policy} $r$, such that
  \begin{equation}
    \label{eq:2a}
    x_{i+1}=r\Bigl(p\bigl(z,\y(\x)\bigr),\x_{1:i},\y_{1:i}\Bigr),
  \end{equation}
  encoding how the algorithm should collect numbers. (Here $\x_{1:i}$ should be understood as the actions taken in the preceding steps $1$ to $i$, and similarly for $\y_{1:i}$). This rule can be simple, for example it could be independent of collected data (regular grids for integration). Or it might have a Markov-type property that the decision at $i$ only depends on $k<i$ previous decisions (for example in ODE solvers). Sometimes, these rules can be shown to be associated with the minimization of some empirical loss function, and thus be given a decision-theoretic motivation. This is for example the case for regular grids in quadrature rules \cite{minka2000deriving}.
\end{enumerate}

\begin{table}
  \centering
  \begin{tabular}{lcccc}
    \toprule
    Problem class & integration & linear opt. & nonlinear opt. & ODE IVPs\\
    \midrule
    inferred $z$ & $z=f; \int f(x)\,dx$ & $z=A^{-1}; Ax=b$ & $z=B=\nabla\nabla\Trans f$ & $z'(t)=f(z(t),t)$\\ 
    classic method & Gaussian quad. & conjugate gradients & BFGS & Runge-Kutta\\
    $p(z)$ & $\GP(f;\mu,k)$ & $\N(A^{-1};M,V)$ & $\GP(z;\mu,k)$ & $\GP(z;\mu,k)$ \\
    $p(\y\g z)$ & $\II(f(x_i)=y_i)$ & $\II(y_i=Ax_i)$ & $\II(y_i = Bx_i)$ & $\II(y_i = z'(t))$\\
    decision rule & minimize & gradient at & gradient under & evaluate at \\
    & post. variance & est. solution & est. Hessian & est. solution \\
    \bottomrule
  \end{tabular}
  \caption{Probabilistic description of several basic numerical problems (shortened notation for brevity). In quadrature, (symmetric positive definite) linear optimization, non-linear optimization, and the solution of ordinary differential equation initial value problems, classic methods can be cast as maximum a-posteriori estimation under Gaussian priors. In each case, the likelihood function is a strict conditioning, because observations are assumed to be noise-free. Because numerical methods are active (they decide which computations to perform), they require a decision rule. This is often ``greedy'': evaluation under the posterior mean estimate. The exception is integration, which is the only area where the estimated solution of the numerical task is not required to construct the next evaluation.}%
\label{tab:models}
\end{table}

The aforementioned results show that classic, base-case algorithms for several fundamental numerical problems can be cast as maximum a-posteriori inference in specific cases of this description; typically under Gaussian priors, and often under simple action rules, like uniform gridding (in quadrature) or greedy extrapolation (in linear algebra, optimization and the solution of ODEs). Table \ref{tab:models} gives a short summary.

\subsection{Current Limitations}

Numerical methods have undergone centuries of development and analysis. The result is a mature set of algorithms that have been ingrained in the scientific tool-set. By contrast, the probabilistic viewpoint suggested here is an emerging area. Many questions remain unanswered, and many aspects of practical importance are missing: Formal analysis is at an early stage. Efficient and stable implementations are still in development. Convincing use-cases from various scientific disciplines are only beginning to emerge. We hope that the reader will take these issues as a motivation to contribute, rather than a hold-up. As the use of large scale computation, simulation, and the use of data permeate the quantitative sciences, there is clearly a need for a formal theory of uncertainty in computation. 

\subsection{New Paths for Research}
\label{sec:new-paths-research}

In our opinion, the match between probabilistic inference and existing numerical methods lays a firm foundation for the analysis of probabilistic numerical methods. We see two primary, complementary goals:

First, implicit prior assumptions can now be questioned. This could be done in an ``aggressive'' way, in the hope of finding either algorithms with faster convergence on a smaller set of problems satisfying stronger assumptions (as in the quadrature example of Section \ref{sec:bayes-herm-quadr}). Conversely, a ``conservative'' re-definition of prior assumptions might improve robustness at increased computational cost. A particularly important aspect in this regard is the action rule $r$. Wherever $r$ is a function of previously collected `data' (known as \emph{adaptive design} in statistics and \emph{active learning} in machine learning), a bias can occur. Where the collected data also influence the result of future actions (as in ODE solvers), a more severe problem, an \emph{exploration-exploitation trade-off} can arise. Checking for such biases, and potentially correcting them, can increase computational cost. But in some applications that require high robustness this effort can pay off.

Secondly the modelling assumptions, in particular the likelihood, can be extended to increase the reach of existing methods to new settings. A first point of interest is the explicit modelling of uncertainty or noise on the evaluations themselves. This generalization, which would be challenging to construct from a classical standpoint, is often straightforward once a probabilistic interpretation has been found. It may be as simple as replacing the point-mass likelihood functions in Table \ref{tab:models} with Gaussian distributions. A prominent case of this aspect is the optimization of noisy functions, such as it arises for example in the training of large-scale machine learning architectures from subsets of a large or infinite dataset. 

Other ideas, like the propagation of knowledge between problems, as in Figure~\ref{fig:CG}, are just as difficult to motivate and study in a classic formulation, but suggest themselves quite naturally in the probabilistic formulation.

There are also practical considerations that shape the research effort. Gaussian distributions play an important role, at least where the inferred quantity is continuous valued. This is not incidental: To a large degree, the point of a numerical method is to turn an intractable computation into a sequence of linear computations. The Gaussian exponential family is closed under linear projections, thus ideally suited for this task.

Efficient adaptation of model hyper-parameters is crucial for a well-calibrated posterior measure. Models with fixed parameters often simply reproduce existing analytic bounds; only through parameter adaptation can uncertainty be actively ``fitted''. Doing so is perhaps more challenging than elsewhere in statistics because numerical methods are ``inner-loop'' algorithms used to solve more complex, higher-level computations. It is important to find computationally lightweight parameter estimation methods, perhaps at the cost of accepting some limitations in model flexibility.

Although the fundamental insight that numerical methods solve inference problems is not new, the study of probabilistic numerical methods is still young. Recent work has made progress, exposing a wealth of enticing applications in the process. We conclude this text by highlighting a most promising, if distant, application motivating ongoing research.

\subsection{A Vision: Chained numerical methods communicating uncertainty}
\label{sec:chain-numer-meth}

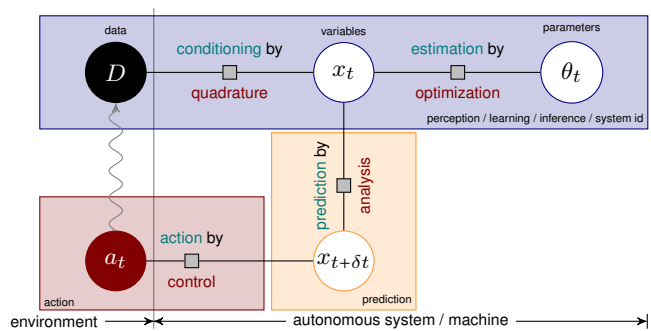
\begin{figure}
  \centering
  {\sf \usetikzlibrary{decorations.text}    
\begin{tikzpicture}

  \draw[gray] (0.5,-3.4) -- (0.5,0.85);

  \draw[<->|] (0.5,-3.3) -- (7,-3.3) 
  node[midway,fill=white] {\scriptsize autonomous system / machine};

  \draw[->] (-1,-3.3) -- (0.5,-3.3)
  node[very near start,fill=white] {\scriptsize environment};


  \fill[fill=dblu,fill opacity=0.2,draw=dblu] (-1,-0.75) rectangle
  (7,0.75);
  \node[anchor=south east,inner sep=0.5mm] at (7,-0.75) {\tiny perception / learning
    / inference / system id};

  \fill[fill=ora,fill opacity=0.2,draw=ora] (2.05,-3.15) rectangle
  (3.95,-0.8);
  \node[anchor=south east,inner sep=0.5mm] at (3.95,-3.15) {\tiny prediction};

  \fill[fill=dred,fill opacity=0.2,draw=dred] (-1,-3.15) rectangle
  (1.95,-1.65);
  \node[anchor=south west,inner sep=0.5mm] at (-1,-3.15) {\tiny action};

  \node[obs] at (0,0) (d) {$D$}; 
  \node[var,draw=dblu] at (3,0) (x) {$x_t$};
  \node[var,draw=dblu] at (6,0) (t) {$\theta_t$};
  \node[var,draw=ora] at (3,-2.5) (c) {$x_{t+\delta t}$};      
  \node[obs,draw=dred,fill=dred] at (0,-2.5) (a) {$a_t$};

  \node[anchor=south] at (d.north) {\tiny data};
  \node[anchor=south] at (x.north) {\tiny variables};
  \node[anchor=south] at (t.north) {\tiny parameters};

  \node[fac] at (1.5,0) (fdx) {} edge (d) edge (x);
  \node[fac] at (4.5,0) (fxt) {} edge (x) edge (t);
  \node[fac] at (3,-1.5)(fxc) {} edge (x) edge (c);
  \node[fac] at (1,-2.5)(fca) {} edge (c) edge (a);

  
  \draw[gray,->,decorate,decoration={coil,aspect=0}] (a) -- (d);

  \node[anchor=south] at (fdx.north) {\scriptsize
    \mpg{conditioning} by};
  \node[anchor=north] at (fdx.south) {\scriptsize
    \dre{quadrature}};

  \node[anchor=south] at (fxt.north) {\scriptsize
    \mpg{estimation} by};
  \node[anchor=north] at (fxt.south) {\scriptsize
    \dre{optimization}};

  \node[anchor=south,rotate=90] at (fxc.west) {\scriptsize
    \mpg{prediction} by};
  \node[anchor=north,rotate=90,text width=2cm,text centered] at (fxc.east) {\scriptsize
    \dre{analysis}};

  \node[anchor=south] at (fca.north) {\scriptsize
    \mpg{action} by};
  \node[anchor=north] at (fca.south) {\scriptsize
    \dre{control}};
\end{tikzpicture}

}
  \caption{\label{fig:pipeline} Sketch of an autonomous system, collecting data to build a parametrised model of the environment. Predictions of future states from the model can be used to choose an action strategy. If intermediate operations are solved by numerical methods, both computational errors and inherent uncertainty should be propagated across the pipeline to monitor and target computational effort.}
\end{figure}

Fueled by ubiquitous collection and communication of data, several academic and industrial fields are now interested in systems that use observations to adapt to, and interact with, their data source in an autonomous way. Figure~\ref{fig:pipeline} shows a conceptual sketch of an autonomous machine aiming to solve a given task by using observations (data) $D$ to build a probabilistic model $p(x_t\g D,\theta_t)$ of variables $x_t$ that describe the environment's dynamics through model parameters $\theta_t$. The model can be used to predict future states $x_{t+\delta t}$ as a function of actions $a_t$ chosen by the machine. The goal is to choose actions that, over time, maximize some measure of utility that encodes the task. This requires a sequence of numerical steps: Inference on $x$ requires marginalization and expectations, i.e.~integration. Fitting $\theta$ involves optimization. Prediction of $x_{t+\delta t}$ may entail solving differential equations. All three areas have linear base cases (inference in linear regression, optimization of quadratic functions, the solution of linear ODEs). The combination of a sequence of ``black-box'' numerical methods in such automated set-ups gives rise to new challenges. Each method receives a point estimate from its precursor, performs its local computation (and adds its local error), and hands the result on. Errors can accumulate in unexpected ways along this chain, but modelling their accumulation provides value: It may be unnecessary to run a numerical method to convergence if its inputs are already known to be only rough estimates.

Specifically, numerical methods allowing for probabilistic inputs and outputs turn the sketch of Figure \ref{fig:pipeline} into a factor graph \cite{Frey97factorgraphs}, and allow propagation of uncertainty estimates along the chain of computation, through message passing \cite{LauritzenSpiegelhalterBP}. This would identify sources of computational error, allowing: the active management of a computational budget across the chain; the dedication of finite computer resources to steps that dominate the overall error; and the truncation of computations early once they reach sufficient precision. Uncertainty propagation through computations has been studied widely before \cite[gives a review]{lee2009comparative}, but the available algorithms focus on the effects of set-up uncertainties on the outcome of a computation, rather than the computation itself. This new functionality explicitly requires calibrated probabilistic uncertainty at each step of the computation, \emph{at runtime}. Classic abstract convergence analyses can not be used for this kind of estimation.

\section{Conclusion}
\label{sec:conclusion}

Numerical tasks can be interpreted as inference problems, giving rise to probabilistic numerical methods. Established algorithms for many tasks can be cast explicitly in this light. Doing so establishes connections between seemingly disparate problems,  yields new functionality, and can improve performance on structured problems. To allow interpretation of the posterior as a statement of uncertainty, care must be taken to ensure well-calibrated priors and models. But even where the uncertainty interpretation is not (yet) rigorously established, the probabilistic formulation already allows for the encoding of prior information about problem structure, including the propagation of collected information among problem instances, leading to improved performance. Many open questions remain for this exciting field. In the long run, probabilistic formulations may allow the propagation of uncertainty through pipelines of computation, and thus the active control of computational effort through hierarchical, modular computations.

\dataccess{N/A. The paper has no data, the empirical parts are all illustrative.}

\conflict{The authors declare to have no conflicting or competing interests.}

\ethics{All three authors contributed content and editing equally.}

\ack{The authors would like to express their gratitude to the two anonymous referees for several helpful comments. They are also grateful to Catherine Powell, University of Manchester, for pointing out the connection to recycled Krylov sequence methods.}

\funding{PH is funded by the Emmy Noether Programme of the German Research Community (DFG). MAG is funded by an Engineering and Physical Sciences Research Council (EPSRC), EP/J016934/2, Established Career Research Fellowship, a Royal Society Wolfson Research Merit Award, and and EPSRC Programme Grant---Enabling Quantification of Uncertainty for Large Scale Inverse Problems---EP/K034154/1.}

\bibliographystyle{rspublicnat}
\bibliography{bibfile}

\end{document}